\documentclass[orivec]{llncs}
\usepackage[latin1]{inputenc}
\usepackage[T1]{fontenc}
\usepackage{graphicx}

\usepackage{algorithm}
\usepackage{algorithmic}

\begin{document}
\title{Solving the Frequency Assignment Problem by Site Availability and Constraint Programming}%
\titlerunning{Site availabity integrated to CP for solving FAP}  

\author{Andréa Carneiro Linhares\inst{1}, Juan-Manuel Torres-Moreno\inst{2,3}, Peter Peinl\inst{4}, Philippe Michelon\inst{2}}%
\authorrunning{Linhares et al.}   
\institute{Universidade Federal do Ceará - Brazil\\
\and
LIA / Universit\'{e} d'Avignon et des Pays de Vaucluse - France\\
\and 
Ecole Polytechnique de Montréal - Canada\\
\and 
University of Applied Sciences Fulda - Germany\\
 {\small andrea.linhares@ufc.br, peter.peinl@informatik.hs-fulda.de, \{juan-manuel.torres, philippe.michelon\}@univ-avignon.fr}
}
\maketitle

\begin{abstract}
The efficient use of bandwidth for radio communications becomes more and more crucial when developing new information technologies and their applications. The core issues are addressed by the so-called Frequency Assignment Problems (FAP).
Our work investigates static FAP, where an attempt is first made to configure a kernel of links.
We study the problem based on the concepts and techniques of Constraint Programming and integrate the site availability concept.
Numerical simulations conducted on scenarios provided by CELAR are very promising. 
\end{abstract}%

\section{Introduction}
\label{intro}

The efficient use of bandwidth for radio communications becomes a more and more crucial aspect in the development of new information technologies and a broad class of its applications. A survey of the literature of the last decade manifests a growing interest in this topic, termed Frequency Assignment Problems (FAP) \cite{AHKMS:03,FPPS:99}. Fundamentally, its objective is to assign frequencies for radio communications avoiding the risk of interference. The dynamic context of FAP was extensively studied in the sense of some combinatorial optimization problems, usually referred to as on-line problems \cite{A:03}. These aim at pointing out the sequential aspect of information arrival and decision making. In the context of FAP, most of the studies related to the on-line algorithms deal with mobile wireless networks \cite{RNW:00}.
Our work studies FAP in radio-based communications networks  \cite{DLAFMV:09,LDFVMD:04}, which is typical of military applications. 
It consists of the allocation of frequencies to radio communications during a tactical military operation. 
Its dynamic dimension, i.e. the progressive establishment of new communication links between existing or new antennas, has been studied in \cite{L:07,D:05}. Using an on-line strategy, frequencies were allocated based on the site availability concept.
We are actually interested in investigating static FAP, issued from the dynamic context. We will study the problem based on the concepts and techniques of Constraint Programming (CP), and integrate the site availability concept. This problem occurs when attempting to configure a kernel, meaning to identify a group of links which are part of the network to be deployed (previously exploited in \cite{D:05}). Certainly, the flexibility of CP is a major asset when modelling the physical constraints of a comunications network. 
In particular, as the constraints in our study are just binary relationships between two frequencies, one may expect rapid responses for a subset of links.
Section \ref{s:etatArt} reviews the state-of-the-art in the domain of frequency assignment applications using CP. 
Following a more precise explanation of the context of our study in section \ref{s:probDescription}, section \ref{s:dispSites} introduces the main idea of the site availability approach and outlines its exact extension. 
The site availability concept is then integrated into the exact solution of our problem in section \ref{s:integrDisp}. 
Our numerical simulations are presented in section \ref{s:ppcTestsBBDisp}. 
Finally, we summarize our conclusions and show some perspectives in section \ref{s:Conclusion}.

\section{State-of-the-art} \label{s:etatArt}

CP is a technique arising from declarative programming and has some roots in computer science and in 
Artificial Intelligence. Recently (around the eighties of last century)
applied to scheduling, planning and assignment problems, CP has attracted the attention of the Operations research community  due to its effectiveness and its potential as a technique to solve combinatorial  optimization problems.
Contrary to Integer Programming, which is based on linear programming, CP mainly favours logic inference.
Inference is realized by the propagation of constraints, which remove infeasible values from domains of variables.
Certainly, the flexibility of this technique is an asset when modelling the relationships among the physical constraints  of a communications network. The characteristics of a constraint are as important as its definition.
We have been interested in pursuing this area of research, becoming more and more important in the literature dealing with the process of solving classical FAPs. 
This choice is also due to the material's properties and the constraints of our problems.
Given the structure of the constraints, i.e. relationships between two frequencies are only binary, one might expect fast results for the assignment of a subset of links. This is due to the observation that the size of kernels in general is small.
A profound vision of CP and its applications has be the subject of a great amount of scientific work. Numerous references are given in \cite{O:00,S:00,M:03}. 

\section{Problem description} \label{s:probDescription}

\begin{figure}[htbp]
  \begin{center}
      \includegraphics[height=5cm]{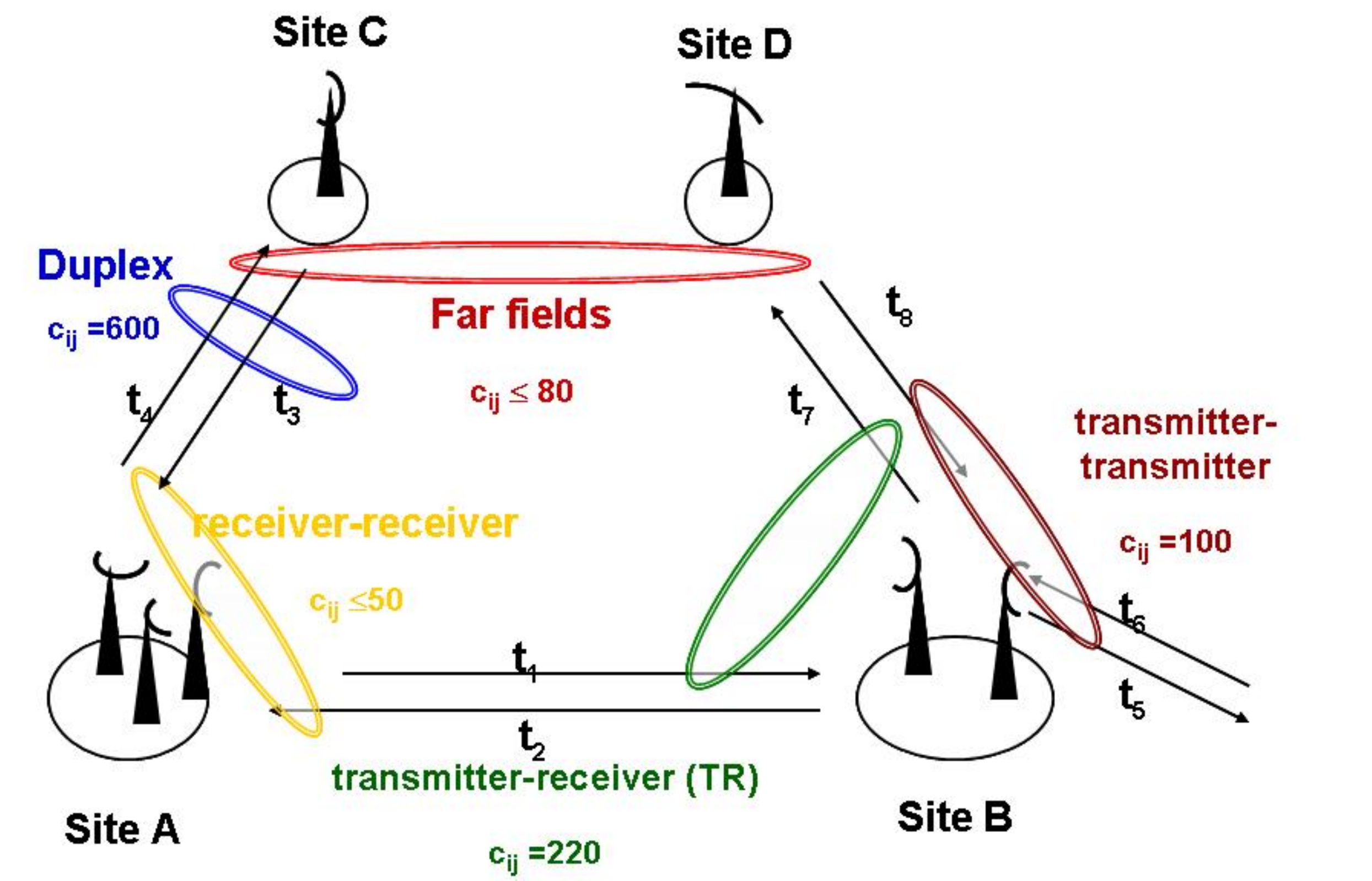}
    \caption[Characteristics of the FAPD constraints] {Characteristics of the FAPD constraints}
    \label{FigCarContraintes}
  \end{center}
\end{figure}

The deployment of a tactical network frequently necessitates the design of a scheme of frequency assignments to satisfy the initial, simultaneous demand of a group of different antennas demanding to establish mutual communication. This group of links is called \textbf{kernel}. The establishment of a link between two sites requires two frequencies (one for each direction of communication). Frequencies are assigned in a way not to cause interference with the group of links that have already been set up. Though all the information about links waiting for frequency assignment is at our disposal the problem conceived is more difficult to solve because of the large number of inherent constraints. Hence, more powerful mechanisms have to be exploited. Despite the strongly combinatorial nature of the problem, there is the same need to come up with an efficient solution in minimal time.
The assignment of frequencies to the kernel (see figure \ref{FigCarContraintes})
is treated in the static context of FAP. Stated more precisely, a pair of sites (origin, destination) will be called a path. A path corresponds to a direction of communication over a link and two sites will never be connected by more than one link. In other words, a link consists of two paths of opposite direction. Let us also point out that the assignment of frequencies to a link is irrevocable, i.e. it will never be questioned in the future.  
We designate the number of sites by $m$, the number of links by $n$, the set of sites by $S=\{s_0,\dots,s_{m-1}\}$, the set of links by $L=\{l_0,\dots,l_{n-1}\}$ and the set of paths by $T=\{t_0,\dots,t_{2n-1}\}$. For each pair of paths $t_i$, $t_j$ in $T$, $c_{ij}$ is the minimum gap between $t_i$ and $t_j$.
$c_{ij}$ is $0$ if no constraint applies to the pair of $t_i$ and $t_j$. For each link $l_i \in L$ , we define $e(t)$ as the transmitter site and $r(t)$ as the receiver one.
The two paths associated with $l_i$ are $t_{2i}$ and $t_{2i+1}$.
Domain D represents the set of frequencies that may be assigned to a link. The values correspond to the french RITA system.
D (identical for all sites) is patitioned into six disjoint sets, called inter-plane, as follows: 
\begin{itemize}
		\item $IPE_0 = {40000,40070,40140}$, $IPE_1 = {41000,41070,41140}$, \\ $IPE_2 = {42000,42070,42140}$, $IPE_3={43000,43070,43140}$,
		\item $IPE_4 = {44000,44070,44140,44210}$, $IPE_5 = {45000,45070,45140,45210}$;
\end{itemize} 

\noindent To avoid noise in the communication (interference), different types of constraints are taken into account. 
They impose a minimum gap between two frequencies assigned to two paths, $t_i$ and $t_{j}$ . 
The difference is represented by $c_{ij}$ and defined by:
$\left|\mathrm f_i-\mathrm f_j\right|\geq \mathrm c_{ij}$, where $\mathrm c_{ij} $ depends on the relationship between $t_i$ and $t_{j}$: 
\begin{itemize}
		\item the duplex constraints between two paths pertaining to the same link impose a gap $c_{ij} = 600$;
		\item the co-site constraints between two paths connected to the same site comprise:
					\begin{itemize}
						\item the transmitter-receiver constraints, where the gap is $c_{ij} = 220$;
						\item the transmitter-transmitter constraints, where the gap is $c_{ij} = 100$;
						\item the receiver-receiver constraints, where the gap $c_{ij}$ rarely exceeds $80$;
					\end{itemize}
		\item the far-field constraints, where the gap $c_{ij}$ practically never exceeds $50$. 
\end{itemize}

\par Hence, the mathematical model of the static problem to be solved is to
find  $f_0,\dots,f_{2n-1}$ such that :
\begin{eqnarray*}
\left| f_i - f_j \right| \geq c_{ij} && (t_i,t_j \in T)\\ f_i \in
D && (t_i \in T)
\end{eqnarray*}

\noindent The maximum number of links of a site is 8. An assignment that fails will be called blockage. With this definition, the objective of our problem is to minimize the amount of blockages in the process of assigning frequency to the kernel.

\section{Site availability} 
\label{s:dispSites}

This strategy concentrates on frequency combinations on the site level. It is based on the realization that, except for the far-field constraints which altogether have little impact (the gap imposed rarely exceeds 50), all constraints concern pairs of paths incident to the same site : co-site, transmitter-transmitter, transmitter-receiver and receiver-receiver. Also, the possibility to assign a frequency to a link connecting two sites $s_1$ and $s_2$ essentially depends on the frequencies assigned to the other incident links ($s_1$, $s_2$). Formally, the availability of a site $s \in S$ for a frequency $f$ is defined as:
\begin{itemize}
	\item the possibility to assign or not assign frequency $f$ to a new transmission path from $s$: transmission availability ($dispE(s,f) \in \left\{ 0,1 \right\} $);
	\item the possibility to assign or not assign frequency $f$ to a new reception path to $s$: reception availability ($dispR(s,f)$ $\in \left\{ 0,1 \right\}$).
\end{itemize}

The idea of the measure defined above is to evaluate the potential intake of the site in its present configuration.
For each new link, that pair of frequencies will be preferred (selected), which conserves maximum availability at the implied sites. A particularity of the context investigated in this paper is that paths always occur in pairs (of the same link). The associated frequencies are designated by $f_{2i}$ and $f_{2i+1}$.
The exact extension of site availability exploits the constraints between all frequencies, available or not, of a site. During the process of allocating frequencies to a new link, new constraints are progressively imposed, for instance the co-site ones, reducing the number of frequencies actually available for subsequent paths to and from the site.
Thus, the exact availability is determined by the maximum number of effective paths in transmission added to the effective maximum number of paths in reception that can be connected to a site. There are two types of exact availability: asynchronous and synchronous. 
For a site $s \in S$, we define asynchronous transmission availability, designated by $dispE^{*}(s)$, as the amount of additional 
transmitters that can be accomodated by the site. Similarly, asynchronous reception availability, $dispR^{*}(s)$, is the amount of additional receivers that can be accomodated by the site. The sum of these two measures gives the exact asynchronous availability of a site, $disp^{*}(s)$.
In the synchronous strategy, we analyze combined reception and transmission availability. We have to find exactly the same number of frequencies in reception and transmission, because we have to choose a pair of frequencies for a new link. Thus, the definition of synchronous availability for a site $s \in S$, designated by $dispC^{*}(s)$, is the amount of additional pairs (receiver, transmitter) that can be held on a site.
This type uses coupling constraints between all available transmitters and receivers (and evaluates all pairs of frequencies). The new measure can be estimated by adding elements to calculate the asynchronous exact availability.

\section{Integration of site availability} 
\label{s:integrDisp}

Site availability is added to search for an exact solution, by employing a decision tree based method. 
This has allowed us to achieve the implantation of mechanisms for the selection of variables, the insertion of an objective into the search for a solution and the definition of additional constraints in the problem model. 
It is well known that a decision tree based method can be improved by modifying the order of the variable and value selection. 
We order the variables by choosing those associated with paths belonging to $Cart8$ sites.
To select the values, we calculate site availability for each frequency in the domain of the path.
In addition, we prohibit those which result in a number of frequencies smaller than the number of paths to assing to the site.
The results obtained in \cite{L:07} suggest that site availability might act as a good lower bound for the solution representing a kernel of links, and yet not increase the complexity of the proposed strategies too much. Thus, we have also used $disp$ as a criterion to evaluate solutions and as a new constraint, in order to filter certain values out of the domain of a path and the sites which it connects. 

\subsection{Pre-processing: reduction of the domain of $Cart8$ sites} 
\label{sb:ppcPreReduCart8}

The existence of sites ($Cart8$) which have to simultaneously react to 8 communication requests makes frequency assignment to a kernel of links more difficult.  We will show that an analysis of the $Cart8$ problem permits to reduce the number of potential solutions.
We know that transmitter-transmitter constraints require a minimum gap of 100 MHz and transmitter-receiver constraints a minimum of 220 MHz. Moreover, paths associated with the same link have to be placed within two differents IPEs to respect a gap of 600 MHz. Receiver-receiver and far-field constraints are very variable. The former impose a gap of at most 80 MHz, but the analysis of our instances (test problems) has shown that real values are about 60 MHz. Taking into account the constraints and the spectrum of frequencies described in section \ref{s:probDescription}, one may conclude that:

\begin{itemize}
	\item Property 1: an IPE will never contain a transmission and a reception path at the same time, because the minimum gap size is 220 MHz;
	\item Property 2: in a small IPE, one may place at most two transmission paths, because the constraint requires a gap of 100 MHz between the two of them; a large IPE, contains two transmission paths at maximum;
	\item Property 3: in a small IPE one may place maximally three receiving paths, in a large IPE up to four receiving paths.
\end{itemize}

\par Let us consider a $Cart8$. For this type of site, we have to place eight transmission and eight reception paths. Conformant to property 2, at least 4 IPEs are required to assign the transmission paths. Therefore, property 1 implies that at most two IPEs may contain reception paths. So, in order to respect property 3, it will be necessary to place the eight reception paths within two IPEs of large size. Subsequently, each small size IPE will get two transmission paths. These properties allow to reduce the path domains in the following way:
\begin{itemize}
\item The domain of a transmission path consists of the frequencies: \\
${40000, 40140, 41000, 41140, 42000, 42140, 43000, 43140}$;
\item The domain of a reception path consists of the frequencies:\\
${44000, 44070, 44140, 44210, 45000, 45070, 45140, 45210}$.
\end{itemize}

\subsection{Algorithm \emph{Branch\&Bound}} 
\label{sb:ppcAlgBBDisp}

The \emph{Branch\&Bound} algorithm constructs the decision tree by first of all evaluating the possibility to find the optimal solution in a particular branch. 
This requires a criterion to evaluate the configurations, which permits to decide whether a partial configuration is more advantageous than another one. This (often) allows to rapidly find a satisfactory solution.
\par Algorithm \ref{algoBB} \cite{D:05} has been used recursively and adapted to find a solution to the kernel problem.
In a second phase, we integrate the site availability approach. The input parameter \emph{nb} corresponds to the number of paths assigned already. Also, by $d$, we designate the domain size, by $t_i$ the pointer to the current path of the search and 
by $D_i[j]$ the j-th value in the domain of path $t_i$.
\emph{reconsider} of a path memorizes the last branch made for each variable.

Function \emph{propager-aff}$(trj,ind)$ propagates the assignment of path $trj$ to the index value $ind$ of its domain. The function returns $true$ if this assignment does not empty any domain of paths in the neighbourhood. In the same way, \emph{propager-des}$(trj)$ propagates the revocation of the assignment of path $trj$. These two functions implement a consistent filtering (remove values from the domains of variables which cannot be part of a solution). This is effected by removing/restoring the values of neighbours that are directly in conflict with new assignments. 

\subsection{Availability as a criterion to evaluate solutions} 
\label{sb:ppcObjBBDisp}

In a first phase, we can use availability as an evaluation criterion of solutions. 
Thus, only those solutions will be retained that leave enough room for a possible dynamic deployment (new links). 

By $best\_disp$ we designate the best avalability value obtained from a branch of the search tree. 
If the value found is better than  $best\_disp$, the current solution is kept and one proceeds with the evaluation according to a predefined ordering of values. 
This test is performed at label $sol$ of algorithm \ref{algoBB}.

\begin{algorithm}[h!]
\caption{Recursive BB()}
\label{algoBB}
\begin{algorithmic}
\STATE {\bf Input} : $nb$  
\STATE {\bf Output} : $s_{sol}$
\STATE {\bf Start}
\STATE {\bf sol}:     
  \IF{$(nb = |L_{trj}|)$}
   \STATE  $s_{sol} \leftarrow$ \emph{save-sol}$()$    
   \STATE  $t_i \leftarrow t_{nb-1}$    
   \STATE  \emph{propager-des}$(t_i)$  
  \ENDIF
  \STATE {\bf next}:  %
  \FOR {$(j = t_i.reconsider +1; j < d; j++)$}
      \IF {($D_i[j]$ is consistent$)$}
          \IF {$($ \emph{propager-aff} $(t_i, j))$}                
		\STATE $t_i.reconsider \leftarrow j$              
		\STATE $nb$++
                 \STATE if {$(nb = |L_{trj}|)$} \emph{goto} {\bf sol}
              \STATE $t_i \leftarrow t_{i+1}$
              \STATE BB()(nb)
              \IF {$(t_i$ \emph{is assigned})}
                     \STATE   \emph{propager-des}$(t_i)$                    
		     \STATE   $nb--$
                     \STATE   \emph{goto} next
              \ENDIF
          \ENDIF
      \ELSE
            \STATE  \emph{propager-des}$(t_i)$       
   \ENDIF
  \ENDFOR
  \STATE {\bf bktk}:
  \STATE $t_i.reconsider \leftarrow 0$   \STATE $t_i \leftarrow t_{i-1}$ \STATE {\bf End}
\end{algorithmic}
\end{algorithm}

\subsection{Availability based selection of values} 
\label{sb:ppcBranchBBDisp}

The selection heuristics determines acceptable values to instantiate variables during the tree search. To integrate availability into this process, immediately on arrival at a new search tree node, availability (transmission and reception) for both sites connected by the path is evaluated. This is done even before assigning the corresponding path.
The calculated availability value is used to evaluate a certain frequency value. While doing so, all values are ruled out that result in a smaller number of frequencies available for assignment than the number of paths of the site that have not yet been
assigned.

\subsection{Availability constraints} 
\label{sb:ppcFiltrageBBDisp}

The integration of availability as a new constraint associated with a site serves as a basis to reduce the domain of variables. This approach distinguishes itself from others, because it works on the domain of sites and takes into consideration the set of paths assigned to each site. 
Let $r$ be the number of paths of a site that have not yet been assigned. Suppose further that transmission and reception availability are calculated for the asynchronous strategy. For each value in the domain, the filtering procedure evaluates $disp^*_E$ and verifies that some frequencies in the neighbourhood of the site are available for assignment. If propagation empties the domain of a neighbour, one removes the frequency from of the path domain.
In addition, if the frequency is not available for any remaining path, it is eliminated from the set of frequencies available for transmission on the site. If, in fact, at least one domain of the neighbours has been emptied by the $r$ assignments, this frequency cannot be allocated to a path on this site. The same process is applied to $disp^*_R$.

\section{Numerical simulations} 
\label{s:ppcTestsBBDisp}
We have decided to perform some simulations on 36 realistic instances (test cases) provided by CELAR ({\sl CEntre d'ELectronique de l'ARmement} (France)). 
The scenarios of these instances contain between 100 and 600 paths (i.e., between 50 and 300 links), between 28 and 168 sites and up to 14755 constraints. Each site has between 1 and 8 links, depending on the scenario. 
The 36 scenarios are divided into 4 groups: instances 01 to 06 (150 links), instances 10 to 19 (50 links), instances 20 to 29 (100 links) and instances 30 to 39 (300 links).

\par The programs have been written in C++. We have performed numerical simulations on a Pentium III double processor (1200 MHz). To validate the strategies outlined in this paper, we have considered exact availability.
We adopt the following notation for our result tables:
\begin{itemize}
\item Availability in the selection of values: $Av-Sel$;
\item Availability in the objective function and in the selection of values: $Av-Obj$;
\item Availability as a constraint, in the objective function and in the selection of values: $Av-Filt$.
\end{itemize}

\par All numerical experiments pre-process $Cart8$ sites. We have used a standard mechanism for the selection of variables. It gives priority to the links which are most restricted. Consequently, paths originating in $Cart8$ sites are put at the top of the list. Filtering turns out to be voracious in terms of processor time and it is only effective when part of the paths has already been assigned. We only employ it for a site with at least 4 links or when half of its paths have already been assigned.

\begin{table} [htbp] 
\centerline{\scriptsize
\begin{tabular}{|c|c|c|c|c|c|c|c|} \hline
 & & \multicolumn{3}{|c|}{$Disp Asynchrone$} & \multicolumn{3}{c|}{$Disp Synchrone$} \\
 & & $PPC-Sel$ & $PPC-Obj$ &  $PPC-Filt$ & $PPC-Sel$ & $PPC-Obj$ & $PPC-Filt$        \\
Instances & $n$ & Max  & Max & Max & Max & Max & Max  \\ \hline
01-06   &   150 & 67.50 (1)& 67.50 (1) & 60.17 (1) & 20.25 (0) & 20.25 (0) & 25.50 (0) \\
10-19   &   50  & 47.20 (8)& 39.30 (4) & 41.90 (6) & 24.80 (2) & 24.80 (2) & 22.85 (2) \\
20-29   &   100 & 41.55 (0)& 40.85 (0) & 40.85 (0) & 28.55 (0) & 28.55 (0) & 30.10 (0) \\
30-39   &   300 & 59.65 (0)& 59.65 (0) & 59.65 (0) & 39.57 (0) & 39.57 (0) & 39.50 (0) \\  \hline
\multicolumn{2}{|c|}{Total}& 53.98     &  51.83    & 50.64     &  28.29    &  28.29    &  29.49 \\  \hline
\end{tabular}}
\caption[Integration of $disp$: resolution for a time of 60s]{Integration of $disp$: resolution for a time of 60s} \label{tab:KMaxPPC60s}
\end{table}

\par Tables \ref{tab:KMaxPPC60s} and \ref{tab:KMaxPPC3600s} compare three different procedures added to the method \emph{Branch\&Bound}, using asynchronous and synchronous availability. 
They show the average number of links assigned
while the solution of our instances is calculated, after computation times of 60 and 3600 seconds.
Incidentally, we identify the number of problems solved exactly, i.e. those in which the number of links assigned corresponds to $n$.
\begin{table} [htbp] \centerline{\scriptsize
\begin{tabular}{|c|c|c|c|c|c|c|c|} \hline
 & & \multicolumn{3}{|c|}{$Disp Asynchrone$} & \multicolumn{3}{c|}{$Disp Synchrone$} \\
 & & $PPC-Sel$ & $PPC-Obj$ &  $PPC-Filt$ & $PPC-Sel$ & $PPC-Obj$ & $PPC-Filt$        \\
Instances & $n$ & Max  & Max & Max & Max & Max & Max  \\ \hline
01-06   &   150 & 67.50 (1) & 67.50 (1) & 67.50 (1) & 40.83 (0) & 40.83 (0) & 26.84 (0) \\
10-19   &   50  & 47.20 (8) & 47.20 (8) & 47.20 (7) & 58.90 (3) & 58.90 (3) & 26.75 (2) \\
20-29   &   100 & 42.30 (0) & 42.30 (0) & 42.30 (0) & 59.70 (0) & 59.70 (0) & 30.90 (0) \\
30-39   &   300 & 61.65 (0) & 61.65 (0) & 61.57 (0) & 79.86 (0) & 79.86 (0) & 40.07 (0) \\ \hline
\multicolumn{2}{|c|}{Total} & 54.66      & 54.66     & 54.64     & 59.82     & 59.82   &  31.14    \\ \hline
\end{tabular}}
\caption[Integration of $disp$: resolution for a time of 3600s]{Integration of $disp$: resolution for a time of 3600s} \label{tab:KMaxPPC3600s}
\end{table}
\noindent The integration of availability permits to solve above all those instances considered to be the most difficult due to the number of $Cart8$ sites, for example instances 20 to 39. 
We cannot compare the results directly to those published in \cite{DLAFMV:09}, because our mechanism of variable selection
makes the problem more difficult in the initial phase. 
Thus, we point out the number of links assigned at 60 seconds. In fact, up to now, there are no numerical simulations similar to those performed in this work on the same 36 scenarios. 
We cannot compare our results to those presented in \cite{D:05}, because the kernel problem solved by our work is not contained in the first 17 links of the scenarios (the links in those scenarios are numbered). 
On the other hand, the solution for the kernel in \cite{D:05}, leans heavily on
repair mechanisms in case the assignment of frequencies should fail (blockage). 
Besides, those scenarios are not public.
They have been provided in the context of the project
{\sl Étude de la Robustesse et d'affectation de fréquences}
and have been made courteously available for other studies and doctoral theses "originating" from the project.

\begin{table} [htbp] 
\centerline{\scriptsize
\begin{tabular}{|c|c|c|c|c|c|c|} \hline
 & \multicolumn{3}{|c|}{$Disp Asynchrone$} & \multicolumn{3}{c|}{$Disp Synchrone$} \\
 & $PPC-Sel$ & $PPC-Obj$ &  $PPC-Filt$ & $PPC-Sel$ & $PPC-Obj$ & $PPC-Filt$        \\
Time & Max    & Max   & Max   & Max   & Max   & Max   \\ \hline
60        & 53.98  & 51.83 & 50.64 & 28.29 & 28.29 & 29,49 \\ 
3600      & 54.66  & 54.66 & 54.64 & 59.82 & 59.82 & 31.14 \\ \hline
\end{tabular}}
\caption{Integration of $disp$: comparaison of global means}
\label{tab:KMaxPPCBilanTemps}
\end{table}

Table \ref{tab:KMaxPPCBilanTemps} shows that algorithm \emph{Branch\&Bound} does not lead to a big improvement, even after 3600 seconds of computation time. 
Besides, the performance of the solution procedures is superior when availability is calculated by the  asynchronous strategy.
This suggests that availability has to remain a guiding principle of the solution algorithm as opposed to a tool to entirely rely on in the search for the exact solution.
Filtering, when used in the synchronous calculation of availability is not efficient. To analyze our experiments with respect to filtering, we have displayed the number of frequency values  filtered. In addition to longer computation times, the mechanism does not achieve to filter out a large number of frequency values. This fact may be directly linked to the pre-processing of $Cart8$ sites, which are taken care of right at the beginning of the search for a solution.

\section{Conclusion and perspectives} 
\label{s:Conclusion}

Using Constraint Programming as a support technique when searching for the solution of the kernel assignment problem introduces important mechanisms into the search for the exact solution. Above all, the idea to base decisions on branching on the study of the physical characteristics of the network permits a higher efficiency throughout the tree search (process). Otherwise the selection of values as well as of variables, turns out to perform better in the analysis of the peculiarities of the instances considered. It seems more interesting to develop specifically adapted procedures or procedures that offer a choice of possibilites when solving a problem, as opposed to rather general procedures.
 We have explored several alternative approaches to obtain the results presented in \cite{L:07}, by integrating site availability for the selection of values when variables are instantiated. This idea has also been exploited for the definition of the objective function and of a new constraint.  Pre-processing instances has achieved a reduction of the frequency possibilities for paths in $Cart8$ sites.
The use of site availability (in a value selection mechanism, in its exact asynchronous version within a \emph{Branch\&Bound} method) shows not to be overly efficient for filtering in the context of our study.
Concerning our strategies, the set of results furnishes some indications to verify. We believe that $ReduPb$ might improve if the IPEs available for a site are retained.
The reduction problem to be solved will be characterized by the assignment of IPEs to link. The advantage of that approach is that two links can be assigned to two large size IPEs of our domain. In practical terms, this observation permits to remove a maximum of 4 links.
That is because a feasible solution can always be easily found.
It also seems interesting, once the existence of a problem without a feasible solution like in the instances 32, 36 et 39, to assign a weight to constraints in the mathematical model. So, the minimization of the number of constraints violated by the solution will become the objective function.
To better evaluate the filtering mechanism that uses availability as a new constraint, we envisage to relax the solution of the exact method, as we have already done when studying greedy approaches. 
When an assignment does not succeed (blockage), the link is removed from the problem and the tree search continues. 
This will permit to determine the number of values filtered, as the problem becomes more "easy" to solve (links minus constraints). 
An interesting perspective might be to include the use of neural networks of the
optimization with constraints \cite{hao} and using an optimal perceptron \cite{torres-moreno:2002}.
This will maximize the objective function by an
evolutionary method dependant on the state of the system.
Data mining algorithms might also be useful in this context.

\bibliographystyle{apalike}
 \bibliography{upao2eng}
\end{document}